\documentclass[11pt,a4paper]{amsart}
\usepackage{amssymb}
\usepackage{amsfonts}
\usepackage{latexsym}
\usepackage{epsfig}
\begin{document}
\newtheorem{thrm}{Theorem}%[chapter]
\newtheorem{thmf}{Th{\'e}or{\`e}me}%[chapter]
\newtheorem{thmm}{Theorem}%[chapter]
\def\thethmm{\ref{twist}$'$}
\newenvironment{thint}[1]{{\flushleft\sc{Th{\'e}or{\`e}me}}
      {#1}. \it}{\medskip} 
\newenvironment{thrm.}{{\flushleft\bf{Theorem}}. \it}{\medskip} 
\newenvironment{thrme}[1]{{\flushleft\sc{Theorem}}
      {#1}. \it}{\medskip} 
\newenvironment{propint}[1]{{\flushleft\sc{Proposition}}
      {#1}. \it}{\medskip} %\def\thethrme{\ref{cls}}
\newenvironment{corint}[1]{{\flushleft\sc{Corrolaire}}
      {#1}. \it}{\medskip} %\def\thethrme{\ref{cls}}
\newtheorem{corr}{Corollary}%[chapter]
\newtheorem{corf}{Corollaire}%[chapter]
\newtheorem{prop}{Proposition}%[chapter]
\newtheorem{defi}{Definition}%[chapter]
\newtheorem{deff}{D{\'e}finition}%[chapter]
\newtheorem{lem}{Lemma}%[chapter]
\newtheorem{lemf}{Lemme}%[chapter]
\def\fy{\varphi}
\def\ul{\underline}
\def\obsf{{\flushleft\bf Remarque. }}
\def\obs{{\flushleft\bf Remark. }}
\def\tl{\tilde}
\def\R{\mathbb{R}}
\def\C{\mathbb{C}}
\def\Z{\mathbb{Z}}
\def\N{\mathbb{N}}
\def\P{\mathbb{P}}
\def\Q{\mathbb{Q}}
\def\p{\pi_1(X)}
\def\Re{\mathrm{Re}}
\def\Im{\mathrm{Im}}
\def\H{\mathbf{H}}
\def\Hs{{Nil}^3}
\def\L{L_\alpha}
\def\h{\frac{1}{2}}
\def\M{\P(E)\times\P(E)^*\smallsetminus\mathcal{F}}
\def\m{\cp2\times{\cp2}^*\smallsetminus\mathcal{F}}
\def\g{\mathfrak{g}}
\def\l{\mathcal{L}}
\def\k{\mathfrak{h}}
\def\s3{\mathfrak{s}^3}
\def\nb{\nabla}
\def\Re{\mathrm{Re}}
\def\Im{\mathrm{Im}}
\newcommand\cp[1]{\mathbb{CP}^{#1}}
\renewcommand\o[1]{\mathcal{O}({#1})}
\renewcommand\d[1]{\partial_{#1}}
\def\sq{\square_X}
\def\gp{\dot\gamma}
\def\j{\mathcal{J}}
\def\iif{\mbox{\bf\em{I}}}
\def\adt{\mbox{\rm ad}_T}
\def\sl2{\mathfrak{sl}_2(\R)}

\title{Automorphism groups of normal $CR$ 3-manifolds}

\author{Florin Alexandru Belgun}
\thanks{AMS classification~: 53A40, 53C25, 53D10.}
\thanks{Supported by the SFB 288 of the DFG}
\date{April 9th, 2001}

\begin{abstract} We classify the normal $CR$ structures on $S^3$ and
  their automorphism groups. Together with \cite{fs}, this closes the
  classification of normal $CR$ structures on contact 3-manifolds. We
  give a criterion to compare 2 normal $CR$ structures, and we show
  that the underlying contact structure is, up to diffeomorphism, unique.
\end{abstract}
\maketitle

\section{Introduction}
A pseudoconvex $CR$ manifold is called {\it normal} if it admits a {\it Reeb}
vector field whose flow preserves the $CR$ structure (in other words,
it admits a symmetry transversal to the underlying contact
structure). If the manifold has dimension 3, this is equivalent to the
fact that it admits a compatible Riemannian {\it Sasakian} structure, for which
the above mentioned Reeb vector field is a unitary Killing vector field. 

In a previous paper \cite{fs} we have classified the Sasakian structures on
compact 3-manifolds and shown that the underlying (therefore normal)
$CR$ structures essentially determine the Sasakian structure (in other
words, the symmetry group of the $CR$ structure is one-dimensional),
in the case where the manifold is not a finite quotient of the
3-sphere. We obtained therefore the classification of the normal $CR$
structures on these manifolds. On the other hand, finite quotients of
$S^3$ admit normal $CR$ structures with bigger automorphism groups
(for example the standard {\it flat} (see next section) $CR$ structure
on $S^3$ admits $PSU(2,1)$ as automorphism group), and in this case
the problem of classifying these structures was still open. 
\smallskip

The first purpose of this paper is to complete the classification of
the normal $CR$ structures on compact 3-manifolds by solving the
remaining case of $S^3$ and its quotients. To do that, we study the
automorphisms of the normal $CR$ structures on $S^3$ and use the 
classification of the Sasakian structures \cite{fs}.

It turns out that the only (locally) homogeneous normal $CR$ structure
is the standard one, and otherwise the (connected component of the)
automorphism group is either a circle or a 2-torus (section 3, Theorem
\ref{main}). There are
other examples of homogeneous $CR$ structures on compact 3-manifolds,
but they are not normal: a left-invariant
$CR$ structure on $S^3$ admits no infinitesimal automorphism {\it
  transversal} to the contact distribution unless it is also
right-invariant, i.e. flat.
\smallskip

% \obs In particular, we retrieve that the automorphism group of a non-flat,
% compact, normal $CR$ 3-manifold is compact --- this holds in general
% for pseudoconvex $CR$ structures

% \obs A direct consequence of that is the following result, similar to
% the well-known theorem of Obata and Lelong-Ferrand in conformal
% geometry: if the automorphism group of a compact $CR$ 3-manifold $M$
% is non-compact, then $M$ is isomorphic to $S^3$ with its canonical
% $CR$ structure. In the general theory of {\it rigid geometric
%   structures} --- e.g. conformal, projective, pseudo-Riemannian, $CR$
% structures --- of Gromov \cite{gro}, \cite{agro}, the above results
% are particular (solved) cases of the {\it vague conjecture}: apart
% from exceptional cases (here: the flat structure), the automorphism
% group of a rigid geometric structure on a compact manifold should be
% compact.

The classification of the normal $CR$ structures on $S^3$ is then
achieved in terms of the {\it wrapping numbers} of a certain uniquely
defined $CR$ Reeb vector field $T_0$, whose orbits are all closed, and
of the metric $g^{T_0}$ induced by projection on the orbit space of
$T_0$ (this space is topologically $S^2$, possibly with 1 or 2 conical
points) --- see section 3 and Theorem \ref{clas}. Note that, {\it a
  priori}, the invariant tensor 
of a $CR$ structure in dimension 3 is the {\it Tanaka curvature
  tensor} (an analogue of the Cotton-York tensor from conformal
geometry), which, even in the simple case of a normal $CR$ structure,
depends on a second order derivative of the curvature of $g^T$.

Finally, we consider the following question: which contact structures
admit compatible normal $CR$ structures? In dimension 3, any contact
manifold $(M,Q)$ admits an infinity of compatible 
$CR$ structures on it~: to pick such a complex structure $J$ on $Q$,
all we need is to choose a positive definite metric on $Q$, as there
is no integrability condition for $J$ required. So,
the existence of a $CR$ structure in dimension 3 implies no
topological condition at all (except orientability), because every
3-manifold admits a contact structure \cite{mar1}. But if we want to
find compatible {\it normal} $CR$ structures on $(M,Q)$, then we
already know that $M$ has to be a
finite quotient of a principal circle bundle over a Riemann surface,
of non-zero Chern class \cite{fs}, \cite{gei}.

Moreover, we want to know what contact structures, up to
diffeomorphism,  on such a manifold admit compatible normal $CR$ structures.
It turns out that there exists at most one, which turns out to be {\it
  tight} (Proposition \ref{cont}).
\smallskip

The paper is organized as follows: in section 2 we recall the
classification of Sasakian structures on $S^3$ and we study in detail
the orbits of the occurring $CR$ Reeb vector fields; in section 3 we
consider, in general, the differential equation related to the
existence of multiple $CR$ Reeb vector fields and describe the
structure of the Lie algebra of the $CR$ automorphism group. The
proof of the main result (Theorem \ref{main}, stated at the beginning
of section 3) is contained in section 4 and 5. The classification
result (Theorem \ref{clas}) follows in section 6 from the previous
sections and the investigation of the flat structures on $S^3$. In
section 7 we describe the contact structures underlying a normal $CR$
structure in dimension 3. 
\smallskip

{\sc Acknowledgments: } The author benefited from useful
discussions with Th. Friedrich and A. Zeghib. 

\section{Preliminaries: Sasakian structures on $S^3$}
In this paper, we call a $CR$ structure $(Q,J)$ on an odd-dimensional
oriented manifold $M^{2n+1}$ with a {\it contact} distribution of hyperplanes
$Q\subset TM$ (such that, if $\eta$ is a non-vanishing {\it (contact)}
1-form for which 
$Q=\ker\eta$, then $\eta\wedge(d\eta)^n$ is a non-vanishing volume
form on $M$), endowed with a {\it compatible} complex structure $J$,
i.e., for $\eta$ as above, $d\eta|_Q=d\eta(J\cdot,J\cdot)|_Q$.

In this case, $d\eta(J\cdot,\cdot)$ is a bilinear symmetric
non-degenerate form on $Q$. If $2n+1=3$, we will always consider
{\it positive} contact forms $\eta$, i.e. such that
$h:=-\frac{1}{2}d\eta(J\cdot,\cdot)$ is a positive definite Hermitian
metric on $Q$ (the renormalization 
factor is useful when considering Sasakian structures --- to be
defined below). The canonical orientation of the manifold is then, by
definition, given by the volume form $\eta\wedge d\eta$.

\obs With this convention, the canonical orientation of $S^3$ with its
standard $CR$ structure induced from its embedding in $\C^2$ is
reverse to the usual one.
\smallskip

In general, in contact geometry one can associate to each contact form
$\eta$ a canonical vector field $V$, called a {\it Reeb} vector field,
such that $\eta(V)=1$ and $\l_V\eta=0$, where $\l$ stands for the Lie
derivative. $V$ is called a $CR$ Reeb vector field if, in addition,
$\l_V|_QJ=0$. Here we follow the convention from \cite{fs} and consider
as a {\it tensor} on a contact manifold any tensorial combination of
elements in the tensorial algebras of $TM$, $Q$, $TM/Q$ and their
duals; with respect to that, $J\in End(Q)$ is a tensor and the Lie derivative
$\l_V$ can be restricted to $End(Q)$ because $Q$ is $V$-invariant, so
the equation $\l_V|_QJ=0$ makes sense. If $V$ satisfies this equation,
it is an infinitesimal automorphism of the $CR$ structure of $M$,
transverse to $Q$; conversely, any such vector field is a $CR$ Reeb
vector field for an appropriate contact form. A Reeb vector field is
{\it positive} if associated to a positive contact form, {\it
  negative} otherwise.

\begin{defi} A {\em CR} structure is called {\em normal} if it admits
  a $CR$ Reeb vector field.
\end{defi}

\begin{defi} A {\em Sasakian} structure $(g,T)$ on a 3-dimensional
  manifold $M$ is a Riemannian metric $g$ admitting a unitary Killing
  vector field $T$ such that $\nabla T$ is a complex structure $J$ on
  $Q:=T^\perp$.
\end{defi}

\obs If $T$ is unitary Killing, then $\nabla T$ is antisymmetric and
vanishes on $T$; the above condition ensures, in addition, that its
square is minus the identity on $Q$. In higher odd dimensions $J$ has
to satisfy some integrability conditions, which are vacuous in
dimension 3.
\smallskip

It is easy to see that, if $(M,g,T)$ is Sasakian, then $(M,Q,J)$, with
the notations above, is a normal $CR$ structure, because $T$ is a $CR$
Reeb vector field. Conversely, it turns out that for any positive $CR$
Reeb vector field $V$, the Riemannian metric defined to be equal to
$h=-\frac{1}{2}d\eta(J\cdot,\cdot)$ on $Q$, for $\eta(V)=1$,
$Q=\ker\eta$, and such that $V$ is unitary and orthogonal to $Q$, is
Sasakian \cite{fs}, \cite{F}.
\medskip

The Sasakian structures on a compact 3-manifold have been completely
classified \cite{fs} using the theory of {\it locally conformally
  K\"ahler metrics} with parallel Lee form on compact complex surfaces
\cite{F}, and, if we could tell which Sasakian metrics on $M$ have the
same underlying $CR$ structure, we could reduce the problem of
classification of the normal $CR$ structures on $M$ to the previous
one, which is solved.

It turns out \cite{fs}, \cite{F}, \cite{gei}, that $M$ is then always
a finite quotient of a non-flat (i.e. with non-zero Chern class)
circle bundle over a Riemann surface (and endowed with the appropriate
orientation), and if this surface has positive genus, the orbits of a
Reeb vector field are always the fibers of such a fibration (The
corresponding Sasakian metrics are called then {\it regular}). In this case, it
follows \cite{fs} that each normal $CR$ structure on $M$ admits, up to
 multiplication by a constant, exactly one $CR$ Reeb vector field,
 and, if we fix the length of its orbits to be equal to $2\pi$, the
 $CR$ structures on $M$ are in 1--1 correspondence to the Sasakian
 structures satisfying the above normalization condition. That closes
 the classification of normal CR structures on 3-manifolds which are
 not covered by $S^3$.

Locally,
 every Sasakian structure can be realized as a fibration over a
 surface, with a suitable connection on it, and we are frequently
 going to use the terminology from fiber bundles: we will often refer
 to the vectors in $Q$ as {\it horizontal}, and to $T$ as {\it vertical}.
\smallskip

From the classification  \cite{fs} we retain the
following facts concerning Sasaki\-an structures on $S^3$ (and its
quotients)~:
On $S^3$ there are regular Sasakian structures (with respect to the
Hopf fibration of $S^3$ over $\cp1$), but also {\it quasi-regular}
(i.e. all the orbits of the Reeb vector field are closed, but of
different lengths) or {\it irregular} (some orbits of the Reeb vector
field are not closed). The structure of the orbits will be discussed
below in detail. 

The orbits of a quasi-regular $CR$ Reeb vector field are the fibers of
a Seifert fibration of $S^3$ over a simply-connected {\it orbifold} with 1
conical point, or with 2 conical points with different angles. In
general, an {\it orbifold} is a topological space endowed with a
differential structure with singularities, and it is given by an atlas
of charts pointing to open sets in $\R^n$ or in quotients of $\R^n$ by
finite linear groups. In dimension 2, a {\it conical} point of {\it
  angle} $2\pi/k$ has a 
neighborhood diffeomorphic to the quotient of the unit disc in $\C$ by
the cyclic group generated by multiplication by a $k$-th root of
unity. It is well-known that the orbifolds that occur above are not
finite quotients of any smooth surface.
\smallskip

{\noindent\bf Convention. } From now on, the manifold denoted by $M$
will be implicitly considered to be diffeomorphic to $S^3$ and
endowed with a normal $CR$ structure. Except otherwise stated, the
symbol $S^3$, when used to denote a $CR$ manifold, will be reserved to the
standard (flat) $CR$ structure on the three-sphere.
\smallskip 

As already mentioned, the Sasakian structures on $M\simeq S^3$ have been
investigated by considering the Riemannian product $M\times S^1$,
endowed naturally with a complex structure (we extend $J$ such that it
sends $T$ into a unitary vector field on the $S^1$ factor);
this complex surface is a {\it primary Hopf} surface {\it without
  resonance (or of class 1)} \cite{F}, i.e. it is the
quotient of 
$\C^2\smallsetminus \{0\}$ by an infinite cyclic group generated by a
contraction 
\begin{equation}\label{hopf}
\gamma(x,y):=(\alpha x,\beta y), \mbox{  where } 0< |\alpha|\le
|\beta|<1,
\end{equation}
and the metric on it is a Hermitian metric with K\"ahler form $\omega$
  and with parallel {\it Lee form} $\theta$ (defined by the
  equation $d\omega=-2\theta\wedge\omega$) 
equal to (a constant times) $ds$, the length form on the factor
$S^1$ \cite{F}. The key information that we will use here is that the
{\it Lee vector field}, which is the dual of the Lee form, and equal
to (a constant times) $-JT$ (where $T$ is the Reeb vector field on the
Sasakian manifold $(M,g,T)$), turns out (\cite{F}, Proposition 8) to
be, up to a constant, identified on $\C^2$ with
$$\log|\alpha|x\partial_x+\log|\beta|y\partial_y,$$
and to get the expression that corresponds to $T$ we have to multiply
this with $i$ and, possibly, a real constant.

The condition to get the above product structure on
$\C^2\smallsetminus \{0\}$ (for example, it is necessary, but not
sufficient, the orbits of $JT$ to be closed) is that $\alpha=|\alpha|$ and
$\beta=|\beta|$, so both are real positive numbers \cite{fs}. The
regularity of the Sasakian structure $(S^3,g,T)$ is then decided by
the ratio $\log\alpha/\log\beta$:
\begin{lem}\label{reg}
$T$ is regular, quasi-regular, resp. irregular if and only if
\linebreak $\log\alpha/\log\beta$ is equal to 1, rational or
irrational. This ratio is 
exactly the ratio of the lengths of the {\em exceptional} orbits of
$T$ contained in the 
complex lines $\{x=0\}$, resp. $\{y=0\}$. If $\log\alpha/\log\beta=p/q$, an
irreducible fraction, the lengths of the orbit in $\{x=0\}$, in
$\{y=0\}$, resp. generic (i.e., different from the 2 exceptional ones)
orbit are $qc$, $pc$ and $pqc$, where $c$ is
a positive constant. If $\log\alpha/\log\beta$ is irrational, then a generic
orbit of $T$ is dense in a 2-dimensional torus; outside the
exceptional orbits, the rest of $S^3$ is foliated by these tori, to
which are tangent all $CR$ Reeb vector fields (which commute with $T$)
and are represented in $\C^2$ by $iax\partial_x+iby\partial_y$, for
$a,b\in\R$ (in particular, this set contains the regular $CR$ Reeb
vector field represented by $ix\partial_x+iy\partial_y$).
\end{lem}
\begin{proof} It is enough to note that the orbits of $T$ in $M$ can
  be identified to the orbits of $i\log\alpha x\partial_x+i\log\beta
  y\partial_y$ on $\C^2$. 
\end{proof}
\obs Let $(M,g,T)$ be a Sasakian structure such that its product
with $S^1$ leads to the Hopf surface $\C^2\smallsetminus
\{0\}/\langle\gamma\rangle$, where
$\gamma(x,y)=\left(a^{1/p}x,a^{1/q}y\right)$, for $p,q\in\N$
mutually prime and $0<a<1$ is a constant. Then $T\equiv - p
x\partial_x-q y\d y$ is 
quasi-regular, and the space $B$ of its orbits in $M\simeq S^3$ (seen as a
factor in the Riemannian product $M\times S^1$, and not as the
standard 3-sphere in $\C^2$) is the space of the complex
orbits of $T$ -- seen as a holomorphic vector field -- in $\C^2\smallsetminus
\{0\}$. It is topologically $\cp1$, and the projection from $\C^2\smallsetminus
\{0\}$ is the holomorphic map $(x,y)\mapsto[x^p:y^q]$, but this fails
to be a submersion at the {\it poles} $[1:0]$ and $[0:1]$. 

To recover the orbifold structure of the basis $B$ we proceed as
follows: around the exceptional orbit contained in $\C\times \{0\}$ we
consider a covering $\C^*\times\C\rightarrow\C^*\times\C$ defined by
$(x,y)\mapsto(x^q,y)$. The pull-back of $T$ is $\tl T=-pqx\d x-qy\d
y$, the space of its orbits $\tl B$ (in $\C^*\times\C$) is $\C$, seen as the
open set in $\cp1$ outside the point $[0:1]$. The corresponding
projection $\C^*\times\C\rightarrow\cp1$ is equivalent with
$(x,y)\mapsto[x^p:y]$ which is a submersion everywhere. The map from
$\tl B$ to $B$, which assigns to each orbit of $\tl T$ the
corresponding orbit of $T$, is then the ramified covering
$[x:y]\mapsto[x^q:y^q]$, which gives the orbifold structure around the
pole $[1:0]$ --- this is thus a conical point of angle $2\pi/q$ (a
smooth point has angle $2\pi$). Note that the orbits of $\tl T$ have
the same length as the generic orbits of $T$, over each of which
project $q$ different orbits of $\tl T$. Over the exceptional orbit
contained in $\C\times\{0\}$ is ``wrapped'' one orbit of $\tl T$ $q$
times. An analogue statement holds for the other pole: $[0:1]$ is
then a conical point of angle $2\pi/p$, there is a ramified covering
with $p$ leaves over it, the source smooth manifold being the orbit
space of $\tl T$, the lift of $T$ on some $p$-fold smooth covering of
$\C\times \C^*$, these orbits generating a fibration, and the generic
orbits of $T$ are covered by $p$ different $\tl T$ orbits, the
exceptional one by only one, ``wrapped'' $p$ times.

The situation where one of $p$ or $q$ is equal to 1 corresponds to the
{\it teardrop} orbifold (a sphere with one conical point), if both are
superior to one we get two conical points.
\smallskip

\begin{defi}The ratios between the lengths of a generic orbit and an
  exceptional one, in case of a quasi-regular $CR$ Reeb vector field $T$,
  are called the {\em wrapping numbers} of $T$.
\end{defi}
Of course, the wrapping numbers coincide with the multiplicity, at the
conical points, of the ramified coverings defined above.

\section{Main statement. The structure of the Lie algebra of
  infinitesimal automorphisms of a normal $CR$ structure on $S^3$}

An important fact, that does not directly follow from the section
above, but from \cite{F}, \cite{fs}, is the following:

% Nevertheless, the {\it local}
% projection $(x,y)\mapsto [x:y^\frac{p}{q}]$, which is well-defined
% (although it depends on some fixed choice of the fractional power)
% everywhere $y\ne 0$, and is a submersion onto a ramified covering of a
% neighborhood of the pole $[1:0]$, with $p$ leaves. The local
% differential (and conformal) model for the pole $[1:0]$ is then the
% quotient of a disk in $\C$ by the multiplication of a primitive $p$-th
% order root of unity, and this gives the orbifold structure (we get an
% analogous model for the other pole) of $B$. The length of an
% orbit of $T$ is proportional with the period of the map $\R\ni
% t\mapsto (xe^\frac{it\pi}{p},ye^\frac{it\pi}{q})$, which is $pq$ if
% both $x,y$ are non-zero, but it is $p$ if $y=0$ and $q$ if $x=0$. If
% we extend the ramified covering map from a neighborhood of $[1:0]\in
% \cp1$ to a neighborhood of the plane $\C^*\times\{0\}$ in
% $\C^2\smallsetminus\{0\}$, it will take the form $(x,y)\mapsto
% (x^p,y^p)$; a generic orbit will be covered by a collection of $p$
% orbits of the same length, while the exceptional orbit will be covered
% by a $p$ times longer orbit
\begin{prop}\label{rapel}
The $CR$  structure of an irregular  Sasakian
metric $(M,g,T)$ always admits a regular $CR$ Reeb vector field $T_0$
which commutes with $T$; on the other hand, if $(g_0,T_0)$ is a regular
Sasakian metric on $S^3$, it admits a $CR$ Reeb vector field commuting
with $T_0$ if and only if the metric induced by $g_0$ on the basis
$\cp1$ (recall that the orbits of $T_0$ are the fibers of a $S^1$
fibration on $\cp1$) admits a Killing vector field.
\end{prop}

To summarize, there are three types of normal $CR$ structures on
$S^3$: with a regular $CR$ Reeb vector field; with a regular and an
irregular $CR$ Reeb vector field, commuting; with a quasi-regular $CR$
Reeb vector field. The last situation can be described, locally around
each Reeb orbit, as (a finite quotient of) a regular Sasakian
structure (see the construction in the remark above; the lifted vector
$\tl T$ is regular). Note that we do not know yet whether these types
overlap or not. In the next sections, we are going to prove:

\begin{thrm}\label{main} If a normal $CR$ structure admits no irregular $CR$ Reeb
  vector field, then it admits a unique (up to a constant) regular or
  quasi-regular $CR$ Reeb vector field; the connected component of its
  $CR$ automorphism group is then $S^1$.

If a normal $CR$ structure admits an irregular $CR$ Reeb vector field,
then it is either the standard flat one (in which case the
automorphism group is $PSU(2,1)$) or it admits an unique regular $CR$
Reeb vector field; in this latter case, the 2 $CR$ Reeb fields span
the (Abelian) Lie algebra of the $CR$ automorphism group, whose
connected component is a 2-dimensional torus.
\end{thrm}

To prove this theorem we will investigate the existence, on a Sasakian
manifold $(S^3,g,T)$, of other $CR$ Reeb vector fields 
\begin{equation}\label{t'}
T'=fT+X_f,
\end{equation}
where $X_f$ is the component in $Q=T^\perp$ of $T'$, and $f$, a
positive function. It turns out \cite{fs} that $T'$ is a Reeb vector field if
and only if 
\begin{equation}\label{xf}
X_f=\frac{1}{2}J(df|_Q)^\sharp,
\end{equation}
and that it is a $CR$ Reeb vector field if and only if 
\begin{equation}\label{ddf}
\sq f:= X.JX.f+JX.X.f-\nb_XJX.f-\nb_{JX}X.f=0,\  \forall X\in Q,
\end{equation}
where $\nabla$ is the {\it Tanaka-Webster} connection on $(M,Q,J,T)$,
defined as follows:
$$\nabla T=0;\ \nabla X\in Q \;\forall X\in Q,\ \nabla J=0,$$
and the following conditions on its torsion $\tau$:
$$\tau(T,X)=0\; \forall X\in Q,\ \tau(X,JX)=2T\|X\|^2,$$
where the norm is relative to the Sasakian metric. The Tanaka-Webster
connection can easily be related to the Levi-Civita connection of the above
Sasakian metric \cite{fs}, but it is more suitable than this one in
the case of a regular 
$CR$ Reeb vector field $T$, because it is easily constructed from the
fibration (following the orbits of $T$) data \cite{fs}; in particular,
the operator $\sq f$, if $f$ is a function constant along $T$, can be
computed using formally the same expression above, but replacing $\nabla$
with $\nabla^B$, the Levi-Civita connection on the basis $B$ (the space of
the orbits of $T$). Such functions appear when we look for $CR$ Reeb
vector fields $T'=fT+X_f$ commuting with $T$, when $T.f=0$ and $X_f$
is the horizontal lift of the vector field $\frac{1}{2}Jdf^\sharp$ on
$B$; the following Lemma holds:
\begin{lem}\label{l1}
The horizontal part $X_f$ of a $CR$ Reeb vector field commuting with
the (quasi)-regular $CR$ Reeb vector field $T$ on $S^3$ projects on the basis
$B$ as a Killing vector field, whose orbits, with the exception of
two, which are points (among which the one or two possible conic
points of the orbifold $B$), are all circles.
\end{lem}
 \begin{proof} The regular case has been treated in \cite{fs}. The
   fact that $X_f$ projects on a Killing vector field follows from a
   local computation, so it holds also in the quasi-regular case. As
   any isometry of an orbifold has to fix (or interchange) the conical
   points, it follows that $X_f$ vanishes at the conical points. Now,
   if a Killing vector field on a (even open) surface has a zero, then it is
   isolated and has index 1; we can smoothen the conical points and
   deform a little $X_f$ such as the resulting vector field on $S^2$
   has the same number of zeros, each with the same index as
   before. They must be thus exactly 2, as the Euler characteristic of
   $S^2$; we will call these degenerate orbits of $X_f$ {\it
   poles}. The other orbits of $X_f$ are smooth submanifolds, because
   they coincide with the level sets of $f$ (recall that
   $X_f=\frac{1}{2}Jdf^\sharp$), and that $df$ is non-degenerate away
   from the 2 poles).
\end{proof}
\begin{lem}\label{irr} Let $T$ be an irregular $CR$ Reeb vector field
  on $S^3$. Then there is a unique regular $CR$ Reeb vector field
  $T_0$ which commutes with $T$.
\end{lem}
\begin{proof} We know that there exists such a $T_0$ \cite{fs}. 
  Suppose $T_0'$ is another regular $CR$ Reeb vector field such that
  $[T,T_0']=0$. Consider the closure of an orbit of $T$: this is a
  torus $\mathcal{T}$. From the standpoint of $T_0'$, $T$ is a
  commuting $CR$ Reeb vector field: as such, $T=f'T_0'+X_f'$, 
  $X_f'$ is a Killing vector field on $B_0'$, the space of the orbits
  of $T_0'$, and we know from Lemma \ref{reg} that the orbits of $T$
  --- as an irregular $CR$ Reeb vector field commuting with $T_0'$ ---
  are dense in tori to which $T_0'$ is tangent. So $T_0'$, as well as
  $T_0$, is tangent to the closures of $T$'s orbits. On a dense set of
  such a torus $\mathcal{T}$, $T_0'$ is determined by its value at a
  point $x$ (as $T_0'$ commutes with $T$, we transport, through the
  flow of $T$, this vector to all
  the points situated on the orbit of $T$ passing through $x$, and
  this set is dense). The
  same holds for $T_0$, so if $T_0'$ coincides with $T_0$ at a point
  of $\mathcal{T}$, it coincides everywhere. In any case we get that
  $[T_0,T_0']=0$ on $\mathcal{T}$, and, as $\mathcal{T}$ was
  arbitrarily chosen, $T_0$ and $T_0'$ commute on $S^3$. As they are
  both regular, they coincide (Lemma \ref{reg}, see also the remark
  below).
\end{proof}

\obs Any $CR$ Reeb vector field has at least 2 closed orbits; for
example, if we consider a $CR$ Reeb vector field $T'$ commuting with a
regular one, $T$, the orbits of $T$ and $T'$ over the two poles $P_1$
and $P_2$ coincide, but they have different lengths; if they both have
length 1 for $T$ (measured using the dual contact form $\eta$), their
lengths for $T'$ (measured using the contact 1-form
$\eta':=f^{-1}\eta$ dual to $T'$) are $f(P_1)^{-1}$,
resp. $f(P_2)^{-1}$. Unless $f$ is everywhere constant, these are then
different, because the values of $f$ at the poles are the only
critical points of $f$, thus its 2 extrema. So there is no other
regular $CR$ Reeb vector field commuting with $T$, and if $T$ admits
other commuting $CR$ Reeb fields, one can always find irregular ones
among them (we can modify the ratio of the maximum and minimum of $f$,
if $f$ is not constant, by adding to $f$ some positive constant). The
last statement also holds if $T$ is a quasi-regular $CR$ Reeb vector
field: any solution $f$ of the equation $\sq f=0, \ X\in Q$ still  has
only 2 critical points (see above), hence they coincide with its
extrema.
\medskip

% If $T$ is 

% The following result
% holds:
% \begin{lem}\label{irr} Let $f$
Of course, the equation $\sq f=0$ being linear, the space of its
solutions is a vector space, containing the constants; this space can
be identified with the Lie algebra of the $CR$ automorphisms group,
the solution corresponding to $T$ is the constant $1$, and the
solutions corresponding to $CR$ Reeb vector fields form the cone of
everywhere-positive solutions. If $T$ is regular, then we can
integrate $f$ on the fibers (the orbits of $T$) using the contact
1-form $\eta$ dual to $T$, and get the function $\iif^{\, T} f$, which is
constant on the orbits of $T$, and still verifies the equation $\sq
\iif^{\, T} f=0$; it thus corresponds to a $CR$ Reeb vector field (more
generally, if $f$ is not everywhere positive, a $CR$ infinitesimal
automorphism) which commutes with $T$. We can consider then
$\iif^{\, T}:\g\rightarrow\g$ as a linear projection on the commutator of
$T$ (it is even a projector if we fix, by no loss of generality, the
length of the orbits of $T$ to be equal to 1). This commutator is
1-dimensional if $B$ admits no Killing vector fields, 2-dimensional if
it admits one, and 4-dimensional if $B\simeq S^2$ has a round metric. 

If $T$ is a quasi-regular $CR$ Reeb vector field, whose generic orbits
have length $pql$, and its exceptional orbits $C_p$ and $C_q$ have
lengths $pl$, resp. $ql$ (for $p,q\in\N^*$ mutually prime, and $l>0$),
and if $f$ is a solution to the equation $\sq f=0, \ X\in Q$, we
define $\iif^{\, T} f$ on a generic orbit as the integral of $f$ on it, and
on an exceptional orbit, say, $C_p$, as $q$ times the integral of $f$
on it. The resulting function $\iif^{\, T} f$ on $B$, the space of the $T$
orbits, is smooth in all the smooth points of $B$, and continuous at
the conical points. Moreover, if we lift it to one of the local
ramified coverings of one of the poles, say, $\{C_p\}$ (see above),
$\tl B_p$, we get a smooth function, which is actually obtained by
integration along the orbits of the regular $CR$ Reeb vector field
$\tl T$ on some smooth covering of a neighborhood of $C_p$ (see
above). It follows that $\iif^{\, T} f$ satisfies the equation for $CR$
infinitesimal automorphisms commuting with $T$, and
$\iif^{\, T}:\g\rightarrow\g$ is, also in this case, a linear projection
(a projector if $pql=1$) on the commutator of $T$, which can be
1-dimensional (no Killing vector field on $B$), or 2-dimensional. As
a non-smooth orbifold cannot admit more than 1 linear independent
Killing fields, these are the only cases.

\obs If a quasi-regular $CR$ Reeb field $T$ admits no linearly
independent commuting $CR$ Reeb field, and if $\sq f=0, \ X\in Q$,
then $\iif^{\, T} f$ is constant everywhere.
\medskip

We are going to use, in the next section, a criterion that will allow
us to conclude, in certain circumstances, that 2 $CR$ Reeb vector
fields $T$ and $T'=fT+X_f$ coincide:

\begin{lem}\label{l2} 1. If $f\equiv 1$ in a neighborhood of an orbit of $T$,
  then $f\equiv 1 $ everywhere.

2. Suppose that the commutator of $T$ in $\g$ is 1-dimensional. Then
   both $T$ and $T'$ are regular or quasi-regular, in the latter case
   with the same exceptional orbits, and same wrapping numbers. If the
   lengths of the generic orbits of both $T$ and $T'$ are equal to 1,
   and if $\iif^{\, T} f=1$ and $\iif^{\, T'} f^{-1}=1$, then $f\equiv 1$ and
   $T\equiv T'$.
\end{lem}
\begin{proof} 
In the first case, $T$ and $T'$ coincide in a neighborhood of one of
their common generic orbits. If they are irregular, we consider the
uniquely (see Lemma \ref{irr}) associated regular $CR$ Reeb vector
fields $T_0$, resp. $T_0'$, and we easily get, as in Lemma \ref{irr},
that they commute in a neighborhood of a torus $\mathcal{T}$, which is
the closure of a common orbit of $T$ and $T'$. Then, for a suitable
constant $k$, $T_0+k[T_0,T'_0]$ is a CR Reeb vector field, and it
coincides with $T_0$ on a neighborhood of an orbit of $T_0$. If we
prove that they coincide everywhere, it follows that $T_0$ and $T'_0$
are two regular CR Reeb vector fields that commute, so (Lemma
\ref{irr} and the Remark thereafter) they coincide.

We can suppose thus that $T$
and $T'$ are both (quasi-) regular, and they coincide on a neighborhood of some
orbit, which we may suppose generic. Then $f$ (the function defining
$T'$ starting from $T$) is 
constant in such a neighborhood, and $\iif^{\,T}f$ is constant on an
open set of $B$, the orbit space of $T$. But $\sq\iif^{\,T}f=0,\forall
X\in TB$, so $J(d\iif^{\,T}f)^\sharp$ is a Killing vector field,
vanishing in an open set, thus everywhere. So $\iif^{\,T}f$  is
constant. At this point we can apply \cite{fs}, Corollary 2, or the
method below, to conclude that $f$ is constant.

The method we apply to prove the second claim in the Lemma is the same
we used --- in the regular case --- in \cite{fs}, Lemma 1; we recall
it briefly: Using the H\"older inequality 
$$\left(\int_{M}f^{-2}\frac{1}{2}\eta\wedge
  d\eta\right)\left(\int_{S^3}f\frac{1}{2}\eta\wedge d\eta\right)^2\ge
\left(\int_{M}\frac{1}{2}\eta\wedge d\eta\right)^3=v^3,$$
we get the implication
$$(\iif^{\,T}f\equiv 1 \mbox{ and } l^T=1)\Longrightarrow v\ge v',$$
with equality if and only if $f$ is constant. Here, $\eta$ is the contact form
such that $\eta(T)\equiv 1$, $v$ is the volume of $M$ with respect
to the Sasakian metric defined by $T$, and $l^T$ is the length of the
orbits of $T$. The hypothesis of the Lemma ensures that the
parenthesis of the above implication, as of the following one, are
true, so we get $v=v'$ and $f=cst.$:
$$(\iif^{\,T'}f^{-1}\equiv 1 \mbox{ and } l^{T'}=1)\Longrightarrow
v'\ge v.$$
If $T$ is quasi-regular, we denote by $l^T$ the length of the generic
orbits of $T$~; besides, all the integrals above will be taken over
the complement of the exceptional orbits in $S^3$. Note that the
exceptional orbits, as well as the wrapping numbers of $T$ and $T'$
coincide (otherwise a linear combination of them would be irregular
--- see the Remark at the end of Lemma \ref{irr} ---,
and this implies that $\g$ is
even-dimensional (see Corollary \ref{dimg} below), but if the
commutator of $T$ is 1-dimensional, then $\g$ is odd-dimensional,
contradiction). The proof follows as in the regular case (we get
equality in the H\"older inequality, thus $f$ is a constant).
\end{proof}
The Lemma above is useful because we can describe Sasakian metrics
\cite{fs}, but have {\it a priori} no simple criterion to check if two
different Sasakian structures are $CR$-isomorphic, or if they admit
other $CR$ Reeb vector fields (except for the commuting ones); as our
goal is to prove, roughly, that there are very few cases when $\dim
G>1$ --- where $G$ is the group of $CR$ automorphisms of $M$ ---, we
will investigate the orbits of $G$ and search for geometric 
facts that would imply the technical hypotheses in Lemma \ref{l2}.
\begin{lem}\label{par}
Let $\g$ be the Lie algebra of infinitesimal automorphisms of a $CR$
structure on $M\simeq S^3$, and let $T\in\g$ be a quasi-regular $CR$ Reeb
vector field. Then the Lie bracket with $T$ is
$\adt:\g\rightarrow\g$, and the integration along the orbits of $T$
yields a projection $\iif^{\,T}:\g\rightarrow \g$. Let $K:=K^T:=\ker\adt$
and $W:=W^T:=\ker\iif^{\,T}$. 

Then $\adt(\g)=\adt(W)=W$ and $\adt|_W$ is an endomorphism of $W$ whose
eigenvalues are all pure imaginary (non-zero). In particular,
$\g=K^T\oplus W^T$, and $W^T$ is even-dimensional.
\end{lem}
\begin{proof} The exponential of $\adt$ is the adjoint action of the
  flow of $T$ on $\g$; as this flow is periodic, it follows that the
  eigenvalues of $\adt$ are imaginary. The image of $\adt$ lies in $W$,
  because, for an infinitesimal $CR$ automorphism $T'=fT+X_f$,
  $\adt(T')=f'T+X_{f'}$, where $f'=T.f$, and, of course, the integral of
  $f'$ along the orbits of $T$ is zero.

It remains to prove that if $T'=fT+X_f\in W\cap K$, then $f\equiv
0$. $T'\in K$ implies that $f$ is constant on the orbits of $T$, and
$T'\in W$ implies that this constant is 0.
\end{proof}
\begin{corr}\label{dimg}
If $(M,Q,J)$ admits an irregular $CR$ Reeb vector field, then
$\dim\g$ is even; if it admits a (quasi-regular) $CR$ Reeb vector
field whose commutator has dimension 1, then $\dim\g$ is odd.
\end{corr}
\begin{proof} It follows from the previous Lemma that $\dim W^T$ is
  always even. The second claim readily follows. For the first claim,
  we take $T$ to be the regular $CR$ Reeb vector field which commutes
 with an irregular one, and then $\dim K^T$ is 2 or 4 (see above).
\end{proof}
 We consider in the
next sections the cases when the commutator of a $CR$ Reeb vector
field is 1 or 2-dimensional (the remaining case, where the commutator
is 4-dimensional, is the flat $CR$ structure on $S^3$; a detailed
study of this structure is contained in Section 5).

\section{Proof of Theorem 1: case with an irregular $CR$ Reeb vector field}

Suppose we have a normal $CR$ structure on $M\simeq S^3$ with an irregular
$CR$ Reeb vector field $T'$; denote by $T$ the (unique) regular one
that commutes to it. We suppose that the $CR$ structure is not flat;
then $\dim K^T=2$ and $\dim\g\ge 2$ is even. We have to prove that
$\dim \g=2$. 
\medskip

If $\dim\g\ge 4$, then all the open orbits of $G$ in $S^3$ are $CR$
flat. Indeed, E. Cartan proved that the automorphism group of a
non-flat 3-dimensional $CR$ manifold has dimension 3 \cite{cart}. If
the union of these orbits is a dense open set in 
$S^3$, the $CR$ structure is flat. The orbits of $G$ already contain
the closures of the orbits of $T'$, which are tori, and 2 exceptional
orbits. If the $CR$ structure is not flat, then, in a neighborhood of
such a torus $\mathcal{T}$, these tori need to coincide with the
orbits of $G$. As $\dim\g>\dim K^T=2$, there is an irregular $CR$ Reeb
vector field $\tl T'$ which does not commute with $T$ or $T'$; denote
by $\tl T$ the corresponding regular $CR$ Reeb vector field, such
that $[\tl T,\tl T']=0$. Let the function $f$ be such that $\tl
T=fT+X_f$. We will prove that, on a neighborhood of $\mathcal{T}$, $f$
is constant on the above mentioned tori. From this it will follow that
$[T,\tl T]$ vanishes on a neighborhood of an orbit of $T$, thus, from
Lemma \ref{l2}, everywhere.

Suppose $f$ is not constant on $\mathcal{T}$. As $\tl T=fT+X_f$ is
tangent to $\mathcal{T}$, it follows that $df(Y)=2g(JX_f,Y)=0$, where
$Y\in Q\cap T\mathcal{T}$; so $f$ is constant on the horizontal curves
that project on circles on $B$, the space of orbits of $T$. If $x$ is
a regular point of $f$ in $\mathcal{T}$, then the horizontal curve
$C_x\subset\mathcal{T}$ containing $x$ is included in a level set of
$f$, thus it is a circle, ``wrapping'' $p$ times over the circle on
$B$ on which $\mathcal{T}$ projects. For neighboring points
$x'\in\mathcal{T}$, the same thing holds, and the length of $C_{x'}$,
in the Sasakian metric determined by $T$, is the same as the length of
$C_x$. If we look at $C_x$ and $C_{x'}$ as the horizontal curves for
the Sasakian metric associated to $\tl T$, they should still have the
same length in this other Sasakian metric; but the horizontal part of
this latter one is multiplied by $f^{-1}$, which is not locally
constant, contradiction.

So $f$ is constant on all the 2-dimensional orbits of  $G$; as we
supposed that their
union contains a non-empty neighborhood of $\mathcal{T}$, then $T$ and
$\tl T$ commute on this open set, thus (see Lemma \ref{l2}, first
point) everywhere. 
\medskip

The Lie algebra $\g$ is thus 2-dimensional and Abelian, so $G$ is a
quotient of $\R^2$. In order to prove that $G$ is actually a torus, we
proceed as follows: let $T$ be the regular $CR$ Reeb vector field,
whose orbits are of length 1, and
let $f$ be a function on the orbit space $B$ of $T$, such that its
differential is the dual of a Killing vector field on $B$, and such
that the image of $f$ is the interval $[0,1]$. Then any infinitesimal
$CR$ automorphism of $M$ is of the 
type $T_{a,b}:=(a+(b-a)f)T+X_{a+(b-a)f}$, where $a$ and $b$ are real
numbers, positive if and only if $T_{a,b}$ is a $CR$ Reeb vector
field. For example $T=T_{1,1}$, and $T_{2,1}$ and $T_{1,2}$ are both
quasi-regular $CR$ Reeb vector fields with one exceptional orbit of
length $1/2$, all the others being of length $1$; the exponentials of
these elements of $\g$ act trivially on $M$, thus they are equal to
$1\in G$. From this it follows that $\exp(T_{p,q})=1\in G,\ \forall
p,q\in\Z$, thus $G\simeq\g/Z^2$ is a torus.

This proves the Theorem in case where there exists an irregular $CR$
Reeb vector field, and the $CR$ structure is not flat.  

\section{Proof of Theorem 1: case without irregular $CR$ Reeb vector fields}

Suppose every $CR$ Reeb vector field on $(M,Q,J)$ is quasi-regular, and
denote by $p:=p^T,q:=q^T$ the wrapping numbers of such a field
$T$. Suppose that the exceptional  
orbits $C_1^T, C_2^T$ have length length $p$, resp. $q$ (if $p$ or
$q$, or if both of them are equal to 1, there is only one, resp. no
exceptional orbit),
and all other orbits have length $pq$, with respect to the dual 1-form
$\eta^T$. It is then clear that, if $p^T$ 
or/and $q^T$ are greater that 1, then the exceptional orbits $C_1^{\tl
  T}$ or/and $C_2^{\tl T}$, and the numbers $p^{\tl T}, q^{\tl T}$ have
to be the same for all $CR$ Reeb vector fields $\tl T$ close to $T$:
this is because the quotient of the lengths of 2 orbits passing
through 2 fixed points $a,b$ is always a rational number, and it must
vary continuously when deforming $T$ to $\tl T$, it is thus
constant; hence $p^{\tl T}=p^T$ and $q^{\tl T}=q^T$ and the orbits of
$\tl T$ passing through points in the 
open set $M\smallsetminus(C_1^T\cup C_2^T)$ have all the same
length, which proves that the exceptional orbits of $\tl T$ coincide
with  the ones of $T$, as claimed.

We study now the orbits of the connected component of $1$ in $G$, the
Lie group of $CR$ automorphisms. The following Lemmas can be proven
using Lemmas \ref{l2} and \ref{par}, by the same method as in \cite{fs}:
 \begin{lem}\label{2orb}
There are no 2-dimensional orbits of $G$.
\end{lem}
% \begin{proof} The ideas used in \cite{fs}, section 4, to prove a
%   similar result on a normal $CR$ manifold which is not a quotient of
%   $S^3$ can be applied also in our case; indeed, the
%   topological restriction in \cite{fs} is used only to ensure that all
%   $CR$ Reeb fields are regular; in our case, this hypothesis is
%   fulfilled in a neighborhood of a hypothetical 2-dimensional orbit
%   (recall that the exceptional orbits are common to all $CR$ Reeb
%   vector fields).
% \end{proof}
\begin{lem}\label{1orb}
There is at most a finite set of 1-dimensional orbits of $G$
\end{lem}
% \begin{proof} Here again, the proof used in \cite{fs}, section 4, can
%   be applied in our case; if we fix a $CR$ Reeb vector field $T$, the
%   union of all these 1-dimensional orbits of $G$ projects on the
%   orbit space $B$ of $T$ in a closed set which, if it is infinite,
%   admits one of its points as an accumulation point of the rest. If
%   this point is a pole of $B$, we can consider, locally, a ramified
%   covering of a neighborhood of the pole, and reduce the problem to
%   the regular case, treated in \cite{fs}.
% \end{proof}

We have thus one open, dense orbit, and a finite number of circular
orbits of $G$, in $M\simeq S^3$. $G$ is odd-dimensional (see Corollary
\ref{dimg}) and, if its dimension is greater that 3, the open dense
orbit (hence the whole manifold) is $CR$ flat \cite{cart}, case which
we exclude (it admits irregular $CR$ Reeb fields). So, if
we suppose that there are at least 2 linearly independent $CR$ Reeb
vector fields, $\dim G=3$, and its Lie algebra is generated by $T$
(that we fix from now on), $V:=fT+X_f$ and $V':=f'T+X_{f'}$, where
$f:M\rightarrow\R$, 
$\iif^{\, T}f=0$ and $f'=T.f$. Moreover, we have $f''=af$, where $a<0$,
thus $f$, restricted to the orbits of $T$, is a sinusoid, and its
critical points are its maximum and minimum. In fact, it has $2pn$
extrema on $C_1$ (of length $p$), $2qn$ extrema on $C_2$, and $2pqn$
extrema on a generic orbit of $T$ ($n\in\N^*$). 
\begin{lem}\label{exorb}
If $\dim G=3$, then $n=1$, and at least one of $p$ or $q$ is equal to
1.
\end{lem}
As the above Lemmas, this follows as in \cite{fs}, from Lemmas
\ref{l2}, \ref{par}, and the following:
\begin{lem}\label{sl2}
If $\dim G=3$, then $\g\simeq \sl2$.
\end{lem}
The proof of that follows as in \cite{fs}, from Lemma \ref{par}, and
the fact that, if $U$ and $V=[T,U]$ span $W^T$, then their bracket
commutes with $T$ and is non-zero, thus is a constant times $T$; we
can then renormalize $T,U,V$ such that they correspond, through an
isomorphism from $\g$ to $\sl2$, to the following matrices:
\begin{equation}\label{tuv}
T\equiv\left(\scriptsize\begin{array}{rr}0&1\\-1&0\end{array}\right)\
U\equiv\left(\scriptsize\begin{array}{rr}1&0\\0&-1\end{array}\right)\
V\equiv\left(\scriptsize\begin{array}{rr}0&1\\1&0\end{array}\right).
\end{equation}

We have thus excluded all but the following situation: $\g\simeq\sl2$,
and the orbits of $G$ are: one circle $C$ and one open orbit $O$, on which we
have a non-flat $CR$ structure. This structure comes then from a
left-invariant $CR$ structure on $G_0:=\tl {SL}_2(\R)$ which extends
smoothly on a compactification (by adding $C$) of the open orbit
$O$. We will prove that this is only possible if the $CR$ structure on
$G_0$ is flat.
\smallskip

% \begin{proof}
% The 3 vector fields $T,V,V'$ must be linearly independent on the
% open orbit of $G$, thus $f$ has no critical point there; on the
% 1-dimensional orbits, $X_f=0$, thus $df|_Q$ vanishes and the critical
% points of $f$ coincide with the extrema on the 1-dimensional orbits 

% \end{proof}
$G$ is a quotient of
$G_0=\tilde{SL}(2,\R)$ by some central subgroup. The centrum of $G_0$
is an infinite cyclic group generated by $e$, and thus
$G=G_p:=G_0/\langle e^p\rangle$, where $p$ is a prime number. For
example, $G_1=PSL(2,\R)$ and $G_2=SL(2,\R)$ (in the latter, $e$
projects on minus the identity). The orbit $O$ is thus isomorphic to a
quotient of $G_0$ by a discrete subgroup $\Gamma$. As $O$ can be
retracted to a circle, $\Gamma\simeq\Z$; as the orbits of $T$ in $O$
are closed, $\Gamma$ must contain a subgroup of the center of $G_0$
(the orbit of $T$ in $G_0$ contains all the center). So $O$ is a
finite quotient of some $G_p,\ p\in\N^*$. % Suppose the quotient is
% taken through a cyclic subgroup generated by a non-central element
% $z$. Then $z\in G$ acts trivially on $O$, hence on the whole $S^3$
\smallskip

{\noindent\bf Important example. } There is a standard action of
$SL(2,\R)$ on $S^3$ which has one open orbit and a 1-dimensional,
circular orbit: consider $S^3$ (with the canonical $CR$ structure) as
embedded in $\C^2$, and let the action of $SL(2,\R)$ be generated by
the following holomorphic vector fields in $\C^2$:
\begin{equation}\label{sl2s3}\begin{array}{ccrcr}
T&\leadsto &-2ix\d x&-&iy\d y\\
U&\leadsto &(1-x^2)\d x&-&xy\d y\\
V&\leadsto &i(1+x^2)\d x&+&ixy\d y
\end{array}\end{equation}
Where $T,U,V$ are the generators of the Lie algebra $\sl2$
(\ref{tuv}), and satisfy the commutation relations 
\begin{equation}\label{comtuv}
[T,U]=-2V,\ [T,V]=2U,\ [U,V]=2T.
\end{equation}
We get these vectors by looking for holomorphic vector fields in
$\C^2$, tangent to $S^3$ (or, equivalently, elements in
$\mathfrak{su}(2,1)$), whose flows preserve the circle $C_0:=\{y=0\}$;
we get a 4-dimensional Lie algebra $\mathfrak{k}$, generated by the above vector
fields and $ix\d x+iy\d y$ (the standard regular $CR$ Reeb vector
field on $S^3$), which commutes with $T$, and then we keep only the
three above, which generate $[ \mathfrak{k},\mathfrak{k}]$; note that
the $CR$ Reeb vector fields contained in this Lie algebra are all
quasi-regular, and all the generic orbits are twice as long as the
unique exceptional orbit $C_0$. For the moment, we only have a Lie
algebra action of $\sl2$ on $S^3$. We easily check that it produces a
Lie group action of $SL(2,\R)$ on $S^3$, and that the local
diffeomorphism $\psi:SL(2,\R)\rightarrow O_0$, where
$O_0:=S^3\smallsetminus C_0$,
defined by $\psi(\mu):=\mu.(0,1)\in O_0,\  \mu\in SL(2,\R)$, is
bijective and has the following expression:
\begin{equation}\label{psi}
\psi(\mu):=\psi\left(\scriptsize\begin{array}{cc}a&b\\
    c&d\end{array}\right) =\left(i\frac{a-d+i(b+c)}{c-b+i(a+d)},
    \frac{2i}{c-b+i(a+d)}\right).
\end{equation}
If we denote by $x_1$, resp. $x_2$, the numerator and the denominator
of the first fraction appearing in the above equation, we get:
$$\psi(\mu)=(x,y)=\left(i\frac{x_1}{x_2},\frac{2i}{x_2}\right),$$
and
$d\psi_\mu=\left(i\frac{dx_1}{x_2}-x\frac{dx_2}{x_2},-y\frac{dx_2}{x_2}\right)$,
from which we conclude that the right-invariant vector fields on
$SL(2,\R)$, whose value in $\mu$ is $T\mu,U\mu,V\mu$, respectively,
are sent through $\psi$ in the vector fields from (\ref{sl2s3}), as
expected, and the left-invariant vector fields equal to $U,V$ in the
identity on $SL(2,\R)$ are sent in:
\begin{equation}\label{drinv}\begin{array}{rrrrr}
\mu U&\leadsto&y^2\d x&-&y^2\overline{\left(\frac{x}{y}\right)}\d y\\
\mu V&\leadsto&iy^2\d x&-&iy^2\overline{\left(\frac{x}{y}\right)}\d y
\end{array}\end{equation}
The fact that the vector field corresponding to $\mu V$ is exactly $i$
times the one corresponding to $\mu U$ is due to the fact that the
left-invariant $CR$ structure $J_1$ induced from $S^3$ on $SL(2,\R)$ by
$\psi$ is such that $U$ and $V$ generate the contact plane $Q$, and
$J_1$ sends $U$ into $V$.
\smallskip

\obs The adjoint action of $T$ on $\sl2$ induces (up to multiplication
by the constant $-1/2$) {\em the same} $CR$ structure on the tangent
space of the identity in $SL(2,\R)$. Actually, $J_1$ is the standard
left-invariant $CR$ structure on $SL(2,\R)=G_2$ (and on its universal
covering $G_0$); any other left-invariant $CR$ structure $J_q, q\ne 0$,
is given --- up to an inner automorphism of the group
--- by the same $Q$ (generated by $\mu U$ and $\mu V$), and such
that 
$$J_q(\mu U):=\frac{1}{q^2}\mu V.$$
Indeed, any element in $\sl2$ whose flow on a (hence on any) $G_p,\
p>0$, is periodic, is equivalent, via an inner automorphism, to $T$;
any plane in $\sl2$ transverse to $T$ and  $\adt$-invariant 
can be brought by an inner automorphism of $\sl2$ to $Q$, and, for the
Sasakian metrics induced by $T$ in the new and in the standard $CR$
structures on $G_0$, the 2 metrics on $Q$ have a common,
left-invariant, orthogonal basis. By an inner automorphism again, we
may bring this basis to $\mu U,\mu V$, as claimed.
\smallskip

We are looking for a left-invariant $CR$ structure on $G$, thus on its
universal covering $G_0$, which, after a quotient by a discrete group,
can be extended smoothly from $O$ to the compact $S^3$. In particular,
if $\tl T$ is a fixed quasi-regular $CR$ Reeb vector field, the metric
on the space of its orbits must have bounded curvature (because this
curvature is determined by the curvature of the Tanaka-Webster
connection corresponding to $T$, which should be well-defined on the
whole compact manifold $S^3$, see \cite{fs}). With no loss of
generality (see the Remark above), we suppose that $\tl T$ is
identified with the right-invariant vector field on $G_2$ $T\mu$,
where $T$ is the matrix (\ref{tuv}) in $\sl2$, and that the
left-invariant $CR$ structure on $G$ is given on $G_2=SL(2,\R)$ by
$(Q,J_q)$ as above. We note that any left-invariant structure on $G$
goes down to any $G_p,\ p>0$, because we take the quotient only by
central elements; it is more convenient to make the computations on
$SL(2,\R)$, because of the example given above; moreover, the space of
the orbits of $T$ in $G_0$ or in any $G_p$, $p>0$, is the same.
\smallskip

First, we know that $\mu U$ and $\mu V$ are orthogonal, that the
ratio of their norms is $1$ in the standard ($CR$ flat) Sasakian
metric associated to $T\mu$ and $J_1$, in general this ratio is $1/q^2$
(for the $CR$ structure $(Q,J_q)$). The product of their norms is half
of the coefficient of $T\mu$ in their commutator, written in the basis
$T\mu,\mu U,\mu V$ (note that they are always linearly independent):
$$[\mu U,\mu V]=-2iy\d y=2\alpha^2 T\mu + ... (\mbox{\it terms in
  Q}),$$
where 
\begin{equation}\label{alph}
\alpha^2:=\frac{2-|y|^2}{|y|^2}.
\end{equation}
We have used the notations from the example above, who realizes
$SL(2,\R)$ as an open set of $S^3\subset \C^2$. In this case, the
projection on the space of orbits of $T\mu$ is
$$\begin{array}{ccccc}
S^3\smallsetminus\{y=0\}&\rightarrow&\cp1\smallsetminus\{[1:0]\}
&\stackrel{\sim}{\rightarrow}&\C\\ 
(x,y)&\mapsto&[x:y^2]&\stackrel{\sim}{\mapsto}&\frac{x}{y^2}.
\end{array}$$
By straightforward computation we get the projections on $\C$ of the
$T\mu $-invariant vector fields $\mu U$ and $\mu V$ to be $\alpha^2\d
z$, resp. $i\alpha^2\d z$, where $z=u+iv$ is the complex parameter on
$\C$ and $\alpha^2:=\sqrt{A}:=\sqrt{1+4|z|^2}$ is the same as before
(\ref{alph}).

For the $CR$ structure given by $J_q$, the vectors
\begin{equation}\label{xy}
X:=X_q:=qf\d u; \ Y:=Y_q:=\frac{1}{q}f\d v,\ f:=\alpha^3=A^\frac{3}{4}
\end{equation}
form an orthonormal frame for the induced metric on the space $\C$ of
the orbits of $T\mu$. We compute
$$[X,Y]=-\frac{1}{q}\d vf X+q\d xfY,$$
then the Levi-Civita connection of the metric for which $X,Y$ is an
orthonormal frame:
$$\begin{array}{ccrcccr}
\nabla_XX&=&\frac{1}{q}\d vfY&\qquad&\nabla_XY&=&-\frac{1}{q}\d vfX\\
\nabla_YX&=&-q\d ufY&\qquad&\nabla_XY&=&{q}\d ufX
\end{array}$$
and its curvature:
\begin{equation}\begin{array}{ccl}
\langle R_{X,Y}Y,X\rangle&=&q^2(f{\d u}^2 f-(\d u
f)^2)+\frac{1}{q^2}(f{\d v}^2 f-(\d v f)^2)\\
&=&q^2(f\Delta f-|df|^2)+\left(\frac{1}{q^2}-q^2\right)(f{\d v}^2f-(\d
vf)^2)\\
&=&q^2\cdot 12A^{-1/2}+\left(\frac{1}{q^2}-q^2\right) 
(6A^{1/2}-48v^2 A^{-1/2}).
\end{array}\end{equation}
The first term in the last line is bounded, the second is not (for
example, on the $Ou$ axis), unless $q=1$, which corresponds to the
standard, flat, $CR$ structure on $SL(2,\R)$.

We have proven that a $\tl{SL}(2,\R)$ action on $O=M\smallsetminus
C\simeq S^3\smallsetminus
\{\mbox{circle}\}$, preserving a normal $CR$ structure on $M\simeq S^3$,
extends on the whole sphere if and only if the $CR$ structure is
flat. As this was the only case where $\dim G$ could be greater that
1, it follows that, if there are no irregular $CR$ Reeb vector fields,
then $\g=\R$ and the connected component of $G$ is a circle. The proof
of Theorem \ref{main} is complete.

\obs Even in the flat case we can classify the above $\tl{SL}(2,\R)$
actions on $S^3$: as we will see in the next section, all regular $CR$
Reeb fields are equivalent (modulo an inner automorphism of $PSU(2,1)$
--- the group of $CR$ automorphisms of the canonical structure of
$S^3$), and all quasi-regular $CR$ Reeb vector fields commute with
exactly one regular one. If we look for the subgroup of $PSU(2,1)$
preserving a circle in $S^3$, we can suppose that this circle is $C_0$
as above, and then get a 4-dimensional Lie algebra, which contains
$\sl2$ as its commutator. We obtain therefore the action described
above, which is an action of $SL(2,\R)$. Note that this is the only
possible action of a group of type $G_p,\, p>0$ prime, satisfying the
required conditions.

\section{$CR$ flat Sasakian structures on $S^3$}
It is known that a $CR$ flat Sasakian structure on a 3-manifold which
is not covered by $S^3$ is --- up to a finite covering --- regular,
and comes from an $S^1$ fibration over a Riemann surface with constant
curvature $k$ \cite{fs}. This is an easy consequence of the fact that,
if a regular Sasakian structure is $CR$ flat, then the {\it Tanaka
  curvature tensor} $\Phi:S^2Q\rightarrow TM/Q$ vanishes. This tensor
is trace-free, and is computed (in the case of a regular Sasakian
structure) from sectional curvature $k$ of the orbit space $B$ of the Reeb
vector field:
$$\Phi(X,X)(T)=-\h\sq k,\; \forall X\in Q,$$
and it vanishes if and only if $J(dk)^\sharp$ is a Killing  vector
field on $B$. If the genus of $B$ is greater than 0, this implies that
$k$ is constant.

If $B$ is a sphere, there are functions $f$ such that $J(df)^\sharp$
is a non-zero Killing  vector field. We are going to investigate the
case where this function is the curvature itself.
\medskip

Let $g$ be a metric on $B\simeq S^2$ with a Killing vector field $Y$, which
vanishes at the {\it poles} $P_1$ and $P_2$. Its orbits are then {\it
  parallels} and we suppose that they have period $2\pi$. The orbits
of $JY$ are then {\it meridians}: geodesics starting from
$P_1$ and pointing towards $P_2$ (we label the poles
as such). Define $X$ to be the unitary vector field on
$S^2\smallsetminus\{P_1,P_2\}$ which is collinear with $Y$ and points
in the same sense. There is a function $r$ such that $Y=rX$. This
function will be viewed as a function on the real line as follows: it
is enough to evaluate it on a meridian, for which we use the arc
length parameterization, the origin of $\R$ corresponds to the pole
$P_1$ (the starting point of the meridian), and $\tau$, the length of
a meridian, corresponds to the pole $P_2$; we extend then
$r:\R\rightarrow\R$ to be an odd function, periodic with period
$2\tau$. The fact that this function is smooth can be easily seen as
follows: we complete a meridian to a closed geodesic; on half of it,
$X$ will be a parallel unitary vector field, and if we extend it to
$\tl X$ to be parallel on the whole closed geodesic (of length
$2\tau$), then $Y(x)=r(x)\tl X(x)$ on half of it, and $Y(x)=-r(x) \tl
X(x)$ on the other half. The function $r:\R\rightarrow\R$ defined
above is just the scalar product between $\tl X$ and $Y$, in the arc
length parameterization of a geodesic. It is smooth.

\obs If the metric on $S^2$ has conical points at the poles (we have
an orbifold -- this happens if the $CR$ Reeb vector field on $S^3$ is
quasi-regular), we can still ``glue together'' $r$ to get a smooth
function: the argument above, applied to a local finite ramified
covering around one pole, identifies that ``glued'' function with a
scalar product restricted to a geodesic. Note that we cannot complete
canonically a meridian in this case, but the function $r$ is constant
on the parallels anyway, so any other such completion, done in a local
ramified covering, yields the desired result.
\smallskip

The fact that the Killing vector field $Y=rX$ is given by the
differential of the curvature $k$ of the metric implies that
\begin{equation}
  \label{eq:dk}
  k'=lr,
\end{equation}
where we have ``glued'' $k$ as above, and considered as a periodic
function on $\R$ (with period $2\tau$ -- the length of a closed
geodesic extending a meridian), the derivative is taken with respect
to the arc length parameterization (i.e. $k'\equiv JX.k$), and $l$ is a
positive real constant (if $l$ were negative, we could change $X$ in
$-X$ and so on).

On the other hand, $k$ can be retrieved, by computation, from $r$: We
use that $\nabla_{JX}JX=\nabla_{JX}X=0$, as $JX$ is unitary and is
tangent to the meridians, which are geodesics, and $[JX,Y]=[JX,rX]=0$,
as $Y$ is Killing and thus preserves the meridians, and get:
$$[X,JX]=\frac{r'}{r}X$$
everywhere except in the poles, thus
$$R_{X,JX}JX=\left(-\frac{r''r-(r')^2}{r^2}-\frac{(r')^2}{r^2}\right)
X= -\frac{r''}{r}X,$$
on the above dense open set in $S^2$. The derivatives are taken with
respect to the arc length parameterization, hence correspond to
derivations along the vector field $JX$. We get thus the following
equation:
\begin{equation}\label{eq:dr}
r''=-kr,
\end{equation}
 which, together with (\ref{eq:dk}), yields a second order
non-linear differential equation on $\R$, for which we seek periodic
solutions:
\begin{equation}\label{eqq}
k''=-\frac{1}{2}k^2+c,\ c\in\R.
\end{equation}
This equation cannot be integrated, in general (using elementary
functions), but it is equivalent to a first order system of ordinary
differential equations, given by the following vector field on $\R^2$
(in coordinates $x,y$):
\begin{equation}\label{Z}
Z:=\left(-\frac{1}{2}x^2+c\right)\d y+y\d x.
\end{equation}
(For any solution $k$ of (\ref{eqq}), the curve $(k,k')$ is an
integral curve of $Z$, and conversely, the first component of any
integral curve of $Z$ is a solution of (\ref{eqq}).)

\obs Any solution to the equation (\ref{eqq}) is a
solution of an equation of the type
\begin{equation}\label{weier}
(z')^2=z^3+pz+q, \ p,q\in \R,
\end{equation}
which comes from the fact that, if $z$ solves (\ref{eqq}), then
$(z,z')$ satisfy an order 3 algebraic equation (see below) that leads
to (\ref{weier}). The general solution to this equation is an elliptic
function and cannot be expressed, in general, in terms of elementary
functions. 

\obs A simple way to find a solution of (\ref{eqq}) could be, {\it a
  priori}, the following: A regular $CR$ Reeb vector field $T$ on
$S^3$ equipped with the standard $CR$ structure is still a Killing
vector field for the Sasakian structure given by another $CR$ Reeb
vector field $\tl T$, provided that $[T,\tl T]=0$; the metrics induced, by
projection, on the space $B$ of orbits of $T$, from the Sasakian
metric on $S^3$ determined by $T$ or by $\tl T$ are different, and their
curvatures $k,\tl k$  both satisfy the equation $\sq \tl k=\sq k=0,\forall
X\in TB$; it is, however, easy to check that $\tl k$ is, as well as $k$,
constant; not only they are conformally equivalent (through the
identity of $B$), but also isometric up to a global homothety. So we
don't get this way any non-trivial solution to (\ref{eqq}).
\medskip

We seek for periodic solutions of (\ref{eqq}), hence for closed
integral curves of $Z$ (\ref{Z}).
First we note that, for negative $c$, the component along $\d y$ of
$Z$ is always negative, fact which excludes any periodic orbit. For
$c=0$ the only periodic orbit is the fixed point $(0,0)$. In general,
the orbits of $Z$ are contained in the level sets of the function
$F:\R^2\rightarrow \R$, 
$$F(x,y):=y^2+\frac{1}{3}x^3-2cx,$$
which are cubics in $\R^2$:

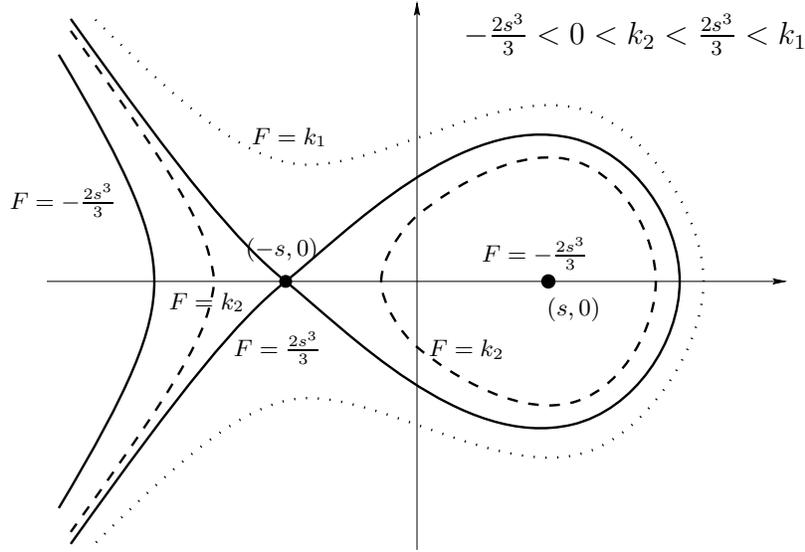
\begin{center}\begin{picture}(0,0)%
\epsfig{file=weierstrass.pstex}%
\end{picture}%
\setlength{\unitlength}{0.00041700in}%
\begingroup\makeatletter\ifx\SetFigFont\undefined%
\gdef\SetFigFont#1#2#3#4#5{%
  \reset@font\fontsize{#1}{#2pt}%
  \fontfamily{#3}\fontseries{#4}\fontshape{#5}%
  \selectfont}%
\fi\endgroup%
\begin{picture}(9762,6924)(2551,-7348)
\put(9376,-4336){\makebox(0,0)[lb]{\smash{\SetFigFont{9}{10.8}{\familydefault}{\mddefault}{\updefault}$(s,0)$}}}
\put(8326,-886){\makebox(0,0)[lb]{\smash{\SetFigFont{12}{14.4}{\familydefault}{\mddefault}{\updefault}$-\frac{2s^3}{3}<0<k_2<\frac{2s^3}{3}<k_1$}}}
\put(5551,-3586){\makebox(0,0)[lb]{\smash{\SetFigFont{9}{10.8}{\familydefault}{\mddefault}{\updefault}$(-s,0)$}}}
\put(5626,-2161){\makebox(0,0)[lb]{\smash{\SetFigFont{9}{10.8}{\familydefault}{\mddefault}{\updefault}$F=k_1$}}}
\put(2551,-2986){\makebox(0,0)[lb]{\smash{\SetFigFont{9}{10.8}{\familydefault}{\mddefault}{\updefault}$F=-\frac{2s^3}{3}$}}}
\put(8551,-3661){\makebox(0,0)[lb]{\smash{\SetFigFont{9}{10.8}{\familydefault}{\mddefault}{\updefault}$F=-\frac{2s^3}{3}$}}}
\put(7876,-4861){\makebox(0,0)[lb]{\smash{\SetFigFont{9}{10.8}{\familydefault}{\mddefault}{\updefault}$F=k_2$}}}
\put(4576,-4261){\makebox(0,0)[lb]{\smash{\SetFigFont{9}{10.8}{\familydefault}{\mddefault}{\updefault}$F=k_2$}}}
\put(5401,-4861){\makebox(0,0)[lb]{\smash{\SetFigFont{9}{10.8}{\familydefault}{\mddefault}{\updefault}$F=\frac{2s^3}{3}$}}}
\end{picture}
\end{center}

For values of $F$ contained between $-\frac{2s^3}{3}$ and
$\frac{2s^3}{3}$, where $s:=\sqrt{2c}$,
its level sets are cubics with 2 connected components in $\R^2$, one
of which is an embedded line, and the other an embedded circle --- each
of which is an orbit of $Z$, because $Z$ does not vanish on these
curves ---; the level set $F=-\frac{2s^3}{3}$ has 2 components, one
of which is the fixed point $(s,0)$ and the other an embedded line ---
which is an orbit of $Z$ ---; the level set $F=\frac{2s^3}{3}$ has
only 1 connected component --- a cubic with a simple double point in
$(-s,0)$ ---, and contains one stationary orbit $\{(0,0)\}$ and 3
different, non-periodic orbits of $Z$, 
which are exactly the connected components of the above cubic after
removal of the fixed point $(-s,0)$ of $Z$. For other values of $F$
its level sets consist, each of them, of one connected cubic in
$\R^2$, diffeomorphic to $\R$, and containing one (non-periodic)
orbit of $Z$. (see figure above)
\smallskip

The only case which is interesting for the equation (\ref{eqq}), for
which we seek periodic solutions, is $c=\frac{1}{2}s^2, \ s>0$ (thus
$c>0$), and the solutions whose orbits coincide with the embedded
circles contained in the level sets of $F$, for values of $F$ between
$ -\frac{2ls^3}{3}$ and $\frac{2s^3}{3}$. Actually all the
solutions obtained this way satisfy the equation (\ref{eqq}), and yield
metrics on $S^2\smallsetminus\{P_1,P_2\}$ with the desired property
($J$ times the differential of the curvature is dual to a Killing
vector field); it remains to check the smoothness at the poles: the
only condition is that the constructed metric gives an angle of
$2\pi$ around each pole. This angle can be measured by integration
around a pole (that is, multiplication by $2\pi$ --- the period of the
Killing vector field $Y=rX$) of the differential of the length of
$Y$. We need thus 
\begin{equation}\label{eq1}
|r'|=1 \mbox{ in } P_1,P_2,
\end{equation}
which is equivalent to
$$\left|-\frac{1}{2}s_i^2+\frac{1}{2}s^2\right|=l,\ i=1,2,$$
where $s_i$ are the solutions, between $-s$ and $s$, of the equation
$F(x,0)=0$, or equivalently, $(s_i,0)$ are the points where the considered
orbit of $Z$ crosses the $Ox$ axis --- these points correspond to
$P_1$ and $P_2$, which are the extrema of $k$, thus where $k'$
vanishes. Because $s_1\in(-s,s)$ and 
$s_2\in(s,2s)$ (the orbit ``turns around'' the fixed point $(s,0)$),
the equation (\ref{eq1}) is equivalent to:
$$s^2-s_1^2=2l \mbox{ and } s_2^2-s^2=2l.$$
On the other hand, the fact that $F(s_1,0)=F(s_2,0)$ (and $s_1\ne
s_2$) is equivalent to
\begin{equation}\label{s12}
s_1^2+s_1s_2+s_2^2=3s^2.
\end{equation}
It follows then that $s_1^2+s_2^2=2s^2$, which leads to
$\sqrt{s^4-4l^2}=s^2$, contradiction.

So there is no metric on $S^2$ such that $J$ times the curvature is
dual to a Killing vector field, except the constant curvature metric
(in which case the corresponding Killing vector field is trivial).
\smallskip

If we are looking for a metric with conical points at the poles
$P_1,P_2$, with wrapping (natural) numbers $q_1$, resp. $q_2$,
the reasoning above remains unchanged, except for the part referring to
the angles given by the metric at the poles; we need to replace
equation (\ref{eq1}) by the following system:
\begin{equation}\label{eq12}
|r'|=\frac{1}{q_1} \mbox{ in } P_1\mbox{ and } |r'|=\frac{1}{q_2}
 \mbox{ in } P_2, 
\end{equation}
which is equivalent to 
$$s^2-s_1^2=\frac{2l}{q_1} \mbox{ and } s_2^2-s^2=\frac{2l}{q_2}.$$
From here, we get the expressions of $s_1,s_2$, then replace them in
(\ref{s12}) and finally get:
$$\left(s^2-\frac{2l}{q_1}\right)\left(s^2+\frac{2l}{q_2}\right)=
\left(s^2-\frac{2l}{q_2}+\frac{2l}{q_1}\right)^2,$$  
which implies
$$s=\sqrt{l\frac{(2/q_2)^3+(2/q_1)^3}{3((2/q_2)^2-(2/q_1)^2)}},\ \mbox{for }
q_1>q_2.$$ 
Thus, for any pair of mutually prime numbers $q_1>q_2$ (the last one
may be equal to 1, in which case we get a teardrop orbifold metric on
$S^2$), and for each fixed $l>0$, there is exactly one value of the
parameter $s$, and exactly 
one orbit of $Z$, thus exactly one solution of (\ref{eqq}), which is
equivalent to
$$k''=\frac{1}{2}(-k^2+s^2);$$
from the point of view of the metric that we obtain on $B\simeq S^2$, the
change in the parameter $l$ means a homothety (i.e., if we fix the
total volume of $B$, then $l$ is unique).

We have proven:
\begin{thrm}\label{plat}
For any regular $CR$ Reeb vector field on a flat $CR$ structure on
$S^3$, the orbit space has constant curvature. The metric on the orbit
space of a quasi-regular $CR$ Reeb vector field $T$ is, up to homothety,
determined by the wrapping numbers of $T$.
\end{thrm}

\obs We obtain, in particular, another proof of the unicity of the
$CR$ flat structure on $S^3$, under the additional supposition that it
admits global $CR$ Reeb vector fields. Of course, this fact holds in
general, as a consequence of the theory of Cartan connections
\cite{tk1}.

\obs As mentioned before, the orbifold metrics on $S^2$ for which the
curvature function $k$ satisfies $\sq k=0,\forall X\in TS^2$ cannot be
expressed in terms of elementary functions {\em in the arc length
  parameterization of the meridians}. In fact, these metrics, an their
curvatures, are easy to compute, as the corresponding Sasakian metrics
on $S^3$, because they are determined by quasi-regular $CR$ Reeb
vector fields on $S^3$, all of which can be explicitely computed. The metrics
on their orbit spaces can be explicitly given in terms of elementary
functions. But this does not lead to any contradiction, because it is the
{\em arc length parameterization of the meridians} --- as a solution
to a second degree differential equation --- which cannot be expressed
in terms of elementary functions. 
% On the other hand, Theorem
% \ref{plat} shows us that it it not necessary to solve explicitly
% equation (\ref{eqq}) or (\ref{weier}) in order to get metrics on (an
% open set of) $S^2$ whose curvature's Hessian is a multiple of the
% identity; we can get such metrics from Sasakian, $CR$-flat, metrics on
% $S^3$, then, by a (local) projection along the orbits of the $CR$ Reeb
% vector field, we get such metrics on (an open set of) $S^2$; it is
% just that we won't get the arc length parameterization of the meridians
% this way. 
\medskip

Theorem \ref{plat} has interesting geometric implications on the
automorphism group of the standard $CR$ structure on $S^3$, namely
$PSU(2,1)$: 
\begin{corr}
The commutator of any regular $CR$ Reeb vector field in\linebreak
$\mathfrak{psu}(2,1)$ is 4-dimensional. For any 2 quasi-regular $CR$
Reeb vector fields $T_1$, $T_2$ with the same wrapping numbers there
is an element in $PSU(2,1)$ mapping $T_1$ on a constant multiple of
$T_2$.
\end{corr}

Another consequence of this is that we can check if two $CR$
structures are isomorphic, just looking at the
$CR$ Reeb vector field and at the metric induced on the orbit space
(recall that any irregular $CR$ Reeb vector field admits exactly one
commuting regular one --- it is easy to find it just by looking at the
lengths of the exceptional orbits of the irregular $CR$ Reeb vector
field --- see section 3), thus reducing the problem of comparing 2
$CR$ structures to the elementary problem of comparing 2 metrics on a
Riemann surface --- more generally, an orbifold:
\begin{thrm}\label{clas}
If 2 Sasakian metrics on a compact 3-manifold $M$ admit the same
underlying $CR$ structure, then:
\begin{enumerate}
\item they are both regular and are {\em 0-type deformation} of each
  other, if $M$ is not a quotient of $S^3$;
\item if $M=S^3$ or a finite quotient, first we replace (if necessary)
  the irregular Sasakian metrics with the (uniquely defined) regular
  Sasakian metrics by {\em basic first type deformations}, so we only need
  to compare (quasi-) regular Sasakian structures. If the Sasakian
  structures are both quasi-regular with the same wrapping numbers
  they are 0-type deformations of each other; If they have different
  wrapping numbers the underlying $CR$ structure is flat, and the
  metric induced on each one's orbit space is determined, following the
  procedure described in the proof of Theorem \ref{plat}, by the
  wrapping numbers. 
\end{enumerate}
\end{thrm} 

By definition \cite{fs}, a {\it 0-type deformation} of a Sasakian
metric is a rescaling of the Reeb vector field (multiplication by a
constant), and a {\it first type deformation} is given by a change in
the $CR$ Reeb vector field (\ref{t'}); it is called {\it basic} if we use in
(\ref{t'}) a function constant along the orbits of the Reeb vector
field. A 0-type deformation is, thus, a particular basic first type
deformation, given by an everywhere constant function. In \cite{F} and
\cite{fs}, a {\it second type deformation} is defined by keeping the
same $CR$ Reeb vector field, the same operator $J$ --- previously
extended to $TM$ by acting trivially on the Reeb vector field ---, and
by changing the metric on the orbit space --- as a consequence, the
contact structure will be modified.
\begin{corr}
A non-trivial second type deformation of a Sasakian structure always
changes the underlying $CR$ structure.
\end{corr}
\begin{proof} If the $CR$ Reeb vector field is quasi-regular, this
  follows from the theorems above; if it is irregular, then any second
  type deformation of it necessarily commutes with the first type
  deformation that replaces the irregular $CR$ Reeb vector field with the unique
  (commuting) regular one, and we can apply the conclusion above. 
 \end{proof}

\section{Contact compact 3-manifolds admitting normal $CR$ structures}
\begin{prop}\label{cont}
  Let $M$ be a compact 3-manifold. Then, up to isomorphism, there
  exists at most one contact structure on $M$ underlying a normal $CR$
  structure. Moreover, this contact structure is {\em tight}.
\end{prop}
By definition, an {\it over twisted} contact structure is one such that
it exists an embedded disc $D\subset M$, transverse to $Q$ everywhere
except in a point $P$ in the interior of $D$, and such that the border
of $D$ is tangent to $Q$. A $CR$ structure is called {\it tight} if it
is not overtwisted. Overtwisted contact structures exist on any
orientable 3-manifold \cite{mar1}, and the isotopy classes of
overtwisted contact structures on $M$ coincide
with their homotopy classes in the category of 2-plane fields in $M$
\cite{e1}. On the other hand, tight structures are more rigid: there
is only one isotopy class of tight contact structures on $S^3$
\cite{e2}; recently, there have been found examples of 3-manifolds not
admitting tight contact structures \cite{eh}.
\begin{proof} We know, from \cite{fs}, that, if $M$ admits a normal
  $CR$ structure, then it is a finite quotient of a circle bundle
  $M_0$ over
  a Riemann surface $B$ with non-zero Chern class, and, in most of the
  cases (see the exceptions below), this $CR$ structure is
  $S^1$-invariant, where the $S^1$-action is free and transverse to
  the underlying contact structure. The latter is thus given by a
  connection on the $S^1$-bundle $M_0\rightarrow B$ (where $B$ is the
  orbit space of $S^1$), such that its curvature is a given volume
  form on $B$ (from here follows that the Chern class should be
  non-zero). On $M_0$, we have thus a $\Gamma$-invariant Sasakian
  structure, where $\Gamma$ is the covering group of the finite
  covering $M_0\rightarrow M$. On $M_0$, all ($\Gamma$-invariant)
  contact structures which are invariant and transverse to the free
  action of $S^1$ are homotopic within ($\Gamma$-invariant) contact
  structures \cite{lutz}, and this ($\Gamma$-equivariant) homotopy is
  followed by a ($\Gamma$-equivariant) isotopy \cite{mart},
  \cite{gray}.

On the other hand, all free $S^1$ actions on $M$ (and, implicitly,
the diffeomorphism type of the quotient space $B$) are isomorphic,
because the only invariants are $b_1(B)$ and the Chern class of the
fibration $M\rightarrow B$, and these are determined by the cohomology
of $M$.
\smallskip

% We can deform it (by changing the
%   metric on $B$ --- on $M_0$ this corresponds to a second type
%   deformation) through $\Gamma$-invariant Sasakian structures, to a
%   {\it standard} one (such that the basis $B$ has constant
%   curvature). The deformation of the underlying contact structures is
%   still $\Gamma$-invariant, and it is followed by a canonically
%   defined isotopy \cite{mart}, which is also $\Gamma$-invariant. This
%   yields a contact isotopy between an arbitrary Sasakian structure on
%   $M$ and a standard one. But all standard Sasakian structures with
%   the same (up to a constant factor) $CR$ Reeb vector field are
%   contact isotopic, and any standard $CR$ Structure on $M_0$ is
%   conjugate to one admitting a fixed $CR$ Reeb vector field.   

The exceptional cases are the normal $CR$ structures on $S^3$ (or
finite quotients) which are $S^1$-invariant, but the $S^1$ action is
not free; it is given by a (unique up to a constant factor)
quasi-regular $CR$ Reeb vector field. We can, nevertheless, deform the
metric on the orbifold $B$ such that it admits a Killing vector field;
the corresponding (second type) deformation of the $CR$ structures is
such that the underlying contact structures are all isotopic
\cite{mart}; but the new $CR$ structure admits at least 2 linearly
independent $CR$ Reeb vector fields, therefore it is invariant to a
free $S^1$ action, which is the case treated above.

On the other hand, the $S^1$ invariant and transversal contact
structures on the total space of a circle fiber bundle as above are
tight \cite{gir}.
\end{proof}

\bigskip

\begin{center}
{\sc {Institut f\"ur Mathematik\\ Humboldt-Universit\"at zu Berlin\\ 
Sitz: Rudower Chaussee 25\\D-10099 Berlin\\Germany}}\\e-mail: {\tt
belgun\@@mathematik.hu-berlin.de}
\end{center}
\bigskip

\end{document}